\documentclass{article}
\usepackage{graphicx}
\usepackage{amsfonts,amssymb,amsmath,amscd, amsthm} 
\title{Alexander Polynomials of Periodic Knots 
\linebreak\Large A Homological Proof and Twisted Extension}
\author{Ross F. Elliot\\
\textit{\small{(rfelliot@gmail.com)}}}
\date{May 5, 2008}

\begin{document}
\newcommand{\Z}[0]{\ensuremath{\mathbb{Z}}}
\newcommand{\T}[0]{\ensuremath{\mathbb{T}}}

\maketitle

\begin{abstract}
\noindent In 1971, Kunio Murasugi proved a necessary condition for a knot to have prime power order.  Namely, its Alexander polynomial $\Delta_K(t)$ must satisfy $\Delta_K(t) \equiv f(t)^{p^r}(1+ t + t^2 + \cdots + t^{\lambda - 1})^{p^r-1}\hspace{5 pt} (mod~p)$ for some knot polynomial $f(t)$ and a positive integer $\lambda$, $(\lambda, p) = 1$.  In this paper I extend this result to the twisted Alexander polynomial.\\

\noindent The direct methods used in the original proof were inadequate for this extension.  Thus, I present an alternate proof using homology theory from which a twisted result follows rather easily.  This result is slightly more complicated than Murasugi's condition, though it has similar features.

\end{abstract}

\section{Introduction}

The following investigation constitutes my undergraduate senior thesis.  It was submitted to Princeton University's Department of Mathematics in partial fulfillment for the degree of Bachelor of Arts.  This effort was supervised by Professor Christopher Skinner and Professor David Gabai who each spent countless hours helping me develop the argument contained herein, as well as teaching me about a wide range of topics only a small fraction of which are used here.  Without this care and concern, I would never have been able to develop my ideas sufficiently to have written a paper of which I am so proud.  I am extraordinarly thankful to both professors for providing me with this formative experience in my mathematical career.  I would also like to acknowledge Hillman, Livingston, and Naik who independently proved this theorem in \cite{Hillman2}.\\

\noindent This content of this paper was initially motivated by the desire to relate ideas in low-dimensional topology and algebraic number theory.  Investigating existing analogies between these areas led to a theorem proved by Kunio Murasugi on periodic knots, which can be thought of as knots that have some sort of rotational symmetry.\\
\newpage
\noindent More precisely, a knot $K$ in $S^3$ is said to be \textit{periodic} of order $n$ if there is an orientation-preserving homemorphism $\phi: S^3\rightarrow S^3$ such that:\\
\indent $1)$ The set of fixed points is a circle (the unknot) disjoint from $K$;\\
\indent $2)$ $\phi(K) = K$;\\
\indent $3)$ $\phi^n = \textbf{1}$ but $\phi^{k} \neq \textbf{1}$ for $0<k<n$.\\ 

\noindent Then for a periodic knot of prime power order, the following condition holds on its Alexander polynomial, a basic invariant of a knot which is computable from a finite presentation of its group \cite{Murasugi2}.\\
\\
\noindent \textbf{Murasugi's Condition.}~~\textit{If $K$ is a periodic knot of order $p^r$ in $S^3$, $p$ a prime, then the Alexander polynomial $\Delta_K(t)$ of $K$ must satisfy: $$\Delta_K(t) \equiv f(t)^{p^r}(1+ t + t^2 + \cdots t^{\lambda - 1})^{p^r-1}\hspace{5 pt} (mod~p)$$
for some knot polynomial $f(t)$ and a positive integer $\lambda$, $(\lambda, p) = 1$.}\\

\noindent This condition relates the cyclotomic polynomial and the Alexander polynomial, two of the most basic objects in algebra and knot theory respectively.  Working with it was far more attractive than with an already disproven conjecture, so we switched gears.  In this paper, I extend Murasugi's condition to the twisted Alexander polynomial, a more complicated knot invariant that is also related to a choice of representation for its group.\\

\noindent \textbf{Twisted Condition.}~~\textit{If $K$ is a periodic knot of order $p^r$ in $S^3$, $p$ a prime, with representation $\rho: \pi_1(S^3-K) \rightarrow GL_n(\Z/p\Z)$ then the twisted Alexander polynomial $\Delta_{K,\rho}(t)$ of $K$ with respect to $\rho$ must satisfy: $$\Delta_{K,\rho}(t) = f(t)^{p^r}\left(\frac{\det(I_n - \rho(l_A)t^\lambda)}{\Delta_{\bar{K}, \bar{\rho}}^0(t)}\right)^{p^r-1}$$ \noindent for some twisted knot polynomial $f(t)$ and a positive integer $\lambda$, $(\lambda, p) = 1$.  Alternatively this condition can be stated as: $$\Delta_{K,\rho}(t) = f(t) \left(\Delta^W_{\bar{K}, \bar{\rho}}(t) \det(I_n - \rho(l_A)t^\lambda)\right)^{p^r-1}$$ \noindent where $\Delta^W$ is another twisted invariant developed by M. Wada.}\\ 

\noindent The relations above are visibly more complicated than those for the regular Alexander polynomial.  The full meaning of the extended condition will be made clear in the course of this paper.\\

\noindent Proving this result was not as simple as following Murasugi's original argument, which became unmanageable when applied to the more complicated twisted case.  Therefore, noticing that the Alexander polynomial of a knot can be defined in terms of the homology of the universal cover of the complement of this knot in $S^3$, I first developed an alternative proof of Murasugi's condition using homology theory.  Then I applied this new argument to the twisted case to obtain my result.  The structure of this paper will be as follows.\\

\noindent Section 2 will recall some useful facts about homology.  Notably, the equivariant homology will provide a means to compute homology from a projective $\Z G$ resolution of $\Z$ given by the free differential calculus, a cornerstone of my argument.  Shapiro's Lemma will be used to relate the homology of a space to that of its cover and will be the justification for an inductive argument.\\

\noindent Section 3 will provide homological definitions for the Alexander and twisted Alexander polynomials, which are usually defined in terms of generators of their elementary ideals in a presentation given by the free differential calculus.  It also includes a brief description of the free calculus and an extension of a theorem of \cite{Hillman}, which given a sequence of modules relates their elementary ideals.\\

\noindent Section 4 presents the homological proof of Murasugi's condition and the following extension to the twisted case.\\

\noindent Limited time prevented me from exploring applications of my result.  An obvious question of interest is whether there exists a knot that satisfies Murasugi's condition but fails to satisfy the condition in the twisted extension for some $q$ (and is thus shown to lack period $q$).  I leave this as an open question to the reader.

\section{Homology}

As both a topological invariant of a space and an algebraic invariant of a group, homology is a convenient tool for studying links between low-dimensional topology and algebra.  Given a topological space $X$, for instance, its first homology group $H_1(X)$ is given by the abelianization of the fundamental group, $\pi_1(X, x_0)$.  Note that the base point $x_0$ is eliminated in the homology since choosing a new base point in a loop will simply permute it cyclically, a distinction that is not maintained in the abelianization.\\

\noindent The truly useful aspect of homology, however, is its formulation in terms of a chain complex $C_*(X)$, which is a sequence of abelian groups connected by homomorphisms.  This has a convenient geometric interpretation.  If $X$ is a simplicial complex ($\Delta$-complex) and $C_n(X) = \Delta_n(X)$ is the free abelian group generated by the $n$-simplices of $X$, then we have the following sequence:

$$\cdots\longrightarrow C_{n+1}\stackrel{\partial_{n+1}}{\longrightarrow} C_{n}\stackrel{\partial_{n}}{\longrightarrow}C_{n-1} \longrightarrow\cdots\longrightarrow C_{1}\stackrel{\partial_{1}}{\longrightarrow} C_{0}\stackrel{\partial_{0}}{\longrightarrow} 0$$

\noindent where $\partial_n: \Delta_n(X) \rightarrow \Delta_{n-1}(X)$ is a boundary homomorphism mapping $n$-simplices into their boundaries, which are $(n-1)$-simplices.  From a direct computation, it is easy to see that $Im\partial_{n+1} \subset Ker\partial_n$, so we can define the $n^{th}$ homology group of $X$ by the quotient $H_n(X) = \frac{Ker\partial_n}{Im\partial_{n+1}}$.\\

\noindent Unless otherwise directed, see \cite{Brown}, \cite{Geoghegan}, \cite{Hatcher} for basic references in this section. 

\subsection{CW Complexes}
Simplicial complexes, however, are not ideal for thinking about a knot $K$ in $S^3$ and its group $G = \pi_1(S^3-K)$.  Therefore, we look to CW complexes which have more algebraic properties.  They are defined inductively as follows.

\newtheorem{prop}{Definition}[subsection]
\begin{prop} \label{xxx}
A CW complex or cell complex is a space $X$ constructed in the following way:\\
\\
\indent $1)$ Let $X^0$ be a discrete set of points, or 0-cells\\
\indent $2)$ If $\left\{e_\alpha^n\right\}_{\alpha\in A}$ is a collection of open $n$-disks, or \textit{$n$-cells}, attach $e_\alpha^n$ to $X^{n-1}$\\
\indent with a map $\varphi_\alpha: S^{n-1} \rightarrow X^{n-1}$ to form $X^n = X^{n-1} \bigsqcup_{\alpha\in A} e_\alpha^n$, which is\\ 
\indent called the \textit{$n$-skeleton} of $X$.\\
\\
If $X = X^n$ for some $n$, then $X$ is finite-dimensional.  The minimum such $n$ is the dimension of $X$.
\end{prop}

\noindent This inductive process is rather intuitive.  Starting with a set of points, form a graph by adding edges and loops.  Then glue open disks onto cycles in the graph.  If two of these open disks are glued to the same cycle, they will form a 2-sphere.  We can then glue open 3-balls to the interiors of these 2-spheres and so forth until we have our $n$-dimensional CW complex.\\

\mbox{}\\
\noindent As defined, CW complexes have an intimate relationship with algebraic structures.  For instance, we can see how the relationships among the lower dimensional skeletons can yield information about the fundamental group of the space.  Geometrically, it is clear that the loop corresponding to a cycle in the 1-skeleton is trivial in the fundamental group if that cycle has a 2-cell glued to it in the 2-skeleton.  In this way, we can view 1-cells as generators and 2-cells as relations so that if we have a presentation of a group $G$ with $s$ generators and $t$ relations, then we can view $G$ as the fundamental group of a space homeomorphic to some CW complex $X$ which has a single 0-cell, $s$ 1-cells, and $t$ 2-cells.  $X$ is then called a \textit{presentation complex} for $G$ \cite{Geoghegan}.  This fact will add understanding to the formulation of the Alexander polynomial in section 3.\\
\mbox{}\\
\noindent Just as simplicial complexes correspond to a theory of simplicial homology, cellular complexes have a corresponding homology theory.  Noting that a quotient $X^n/X^{n-1}$ of skeletons corresponds to the $n$-cells of a CW-complex $X$ we obtain the following definition.

\begin{prop} \label{xxx}
The complex: $$ \cdots \longrightarrow H_{n+1}(X^{n+1}, X^n) \stackrel{\partial_{n+1}}{\longrightarrow} H_n( X^n, X^{n-1} ) \stackrel{\partial_n}{\longrightarrow} H_{n-1}( X^{n-1}, X^{n-2} ) \longrightarrow \cdots$$ is called the cellular chain complex of $X$.  The cellular homology of $X$ is then the homology of its cellular chain complex.
\end{prop}

\subsection{Equivariant Homology}

There is a more algebraic extension of the homology theories discussed so far which gives a generalized module structure to the chain complex.  This structure inherently allows certain actions on the space whose importance will become apparent in section 4.\\

\noindent A preliminary to understanding this homology with coefficients and its specialization to equivariant homology is the concept of a projective module.  A module $P$ is said to be \textit{projective} if for modules $M$, $M'$ and for every homomorphism $\varphi : P \rightarrow M$ and every surjective homomorphism $i: M' \rightarrow M$ there is a homomorphism (lift) $\psi: P \rightarrow M'$ such that $i\psi = \varphi$.  We can now define a projective resolution which is related intimately with the aforementioned homology theories.

\begin{prop} \label{xxx}
Let $R$ be a ring and $M$ be a left R-module.  A projective resolution of M is an exact sequence of $R$-modules $$\cdots\longrightarrow P_2\stackrel{\partial_2}{\longrightarrow} P_1\stackrel{\partial_1}{\longrightarrow} P_0 \stackrel{\epsilon}{\longrightarrow} M\longrightarrow 0$$ such that each $P_i$ is a projective module or, equivalently, is the direct summand of a free module.
\end{prop}

\noindent The homology with coefficients is the homology of a complex formed by tensoring a projective resolution with a module.  If this module is $\Z$, then the homology with coefficients agrees with the regular homology, so the extension is natural.  This process is defined as follows.

\begin{prop} \label{hwc}
For a group $G$ let $M$ be a $G$-module and $F$ be a projective resolution of $\Z$ over the group ring $\Z G$.  Then $H_*(G; M)$, the homology of $G$ with coefficients in $M$, is given by $H_*(F\otimes_G M)$.
\end{prop}

\noindent The equivariant homology is then a specialization to the case $F = C_*(X)$ where $C_*(X)$ is the chain complex of a CW-complex $X$ with an associated $G$-action on $X$ which freely permutes its cells.  This type of space is called a \textit{$G$-complex.}

\begin{prop} \label{xxx}
Let $C(X)$ be the cellular chain complex of a $G$-complex $X$.  Then $H^G_*(X)$, the equivariant homology of $X$, is given by $H_*(G; C(X))$ \cite{Brown}.  Furthermore, if $M$ is a $G$-module, then there is a diagonal $G$-action on $C(X) \otimes_{\Z G} M$, and $H^G_*(X; M)$, the equivariant homology of $X$ with coefficients in $M$, is given by $H_*(G; C(X) \otimes_{\Z G} M)$.   
\end{prop}

\noindent A finite, connected CW-complex $X$ with $G = \pi_1(X)$ will not, in general, have the structure of a $G$-complex.  However, its universal cover $\tilde{X}$ will have the structure of a $G$-complex (so $C_*(\tilde{X})$ is a projective $\Z G$-resolution of $\Z$).  Henceforth, we define the notation: $$H_*(X; M) := H_*^G(\tilde{X}; M).$$  

\noindent To calculate the groups $H_*(X; M)$, we can use any projective $\Z G$-resolution $F$ of $\Z$, not just $C_*(\tilde{X})$.  This fact will be useful later on when we use a resolution coming from the free differential calculus.

\subsection{Shapiro's Lemma}

Originally proved by Arnold Shapiro at the request of Andre Weil, Shapiro's Lemma relates the coefficient homology of a group to that of a subgroup \cite{Brown}, \cite{Geoghegan}.

\newtheorem{lem}{Lemma}[subsection]
\begin{lem}[Shapiro] \label{xxv}
Let $H \subseteq G$ and let $M$ be an $H$-module.  Then $$H_n(H; M)\cong H_n(G; \Z G \otimes_{\Z H} M)$$
\end{lem}

\begin{flushleft}
There is a topological analog of the lemma that relates the equivariant homology of a space to that of its cover.
\end{flushleft}

\newtheorem{cor}{Corollary}[subsection]
\begin{cor} \label{xxv}
Let $X$ be a space with $G = \pi_1(X)$.  Corresponding to a subgroup $H \subset G$ is a cover $\bar{X}$ of $X$ such that $\pi_1(\bar{X}) = H$.  Further, let $\tilde{X}$ be the universal cover of $X$ (and therefore $\bar{X}$), and let $M$ be a $G$-module.  Then $H_*(\bar{X}; M)$ is the homology of $C_*(\tilde{X}) \otimes_{\Z H} M$, which is given by $$H_n(\bar{X}; M)\cong H_n(X, Ind_H^G M)$$ where $Ind_H^G M = \Z G \otimes_{\Z H} M \cong \Z [G/H] \otimes_\Z M$.
\end{cor}

\noindent This latter form will be useful for the proof in section 4.  Henceforth, by "Shapiro's Lemma," I will mean the topological form.

\subsection{Mayer-Vietoris Sequences}

Often the direct calculation of specific homology groups is tedious or infeasible.  In this case, it is convenient to have a method by which one can express them in terms of the homology of spaces with known or more easily computable homology groups.  A \textit{Mayer-Vietoris sequence} does just this by associating a decomposition of a space into two subspaces with a long exact sequence of homology groups.\\

\noindent To derive this sequence, first let a space $X$ be the union of the interiors of some subspaces $U$ and $V$.  Also let $C_*(U + V)$ be the subgroup of $C_*(X)$ composed of sums of chains in $U$ and chains in $V$.  Then the usual boundary maps on $C_*(X)$ are also boundary maps on $C_*(U + V)$ and thus, the latter is a chain complex.  Letting $\phi(x) = (x, -x)$ and $\psi(u, v) = u + v$, we have that: $$0\longrightarrow C_*(U \cap V) \stackrel{\phi}{\longrightarrow} C_*(U) \oplus C_*(V) \stackrel{\psi}{\longrightarrow} C_*(U+V) \longrightarrow 0$$ is a short exact sequence.  Finally, by Proposition 2.21 of \cite{Hatcher}, the inclusion $C_*(U+V)\hookrightarrow C_*(X)$ is a chain homotopy equivalence that induces an isomorphism $H_*(U+V)\cong H_*(X)$.  Thus, we have the following definition.

\begin{prop} \label{xxx}
Let $U$, $V$ be two open subspaces of a space $X$ such that $X=U \cup V$. Then the Mayer-Vietoris sequence associated to this decomposition is the long exact sequence: $$\cdots \rightarrow H_n(U \cap V) \rightarrow H_n(U) \oplus H_n(V) \rightarrow H_n(X) \rightarrow H_{n-1}(U \cap V)\rightarrow\cdots\rightarrow 0$$ which is obtained from the aforementioned short exact sequence of chain complexes.  
\end{prop}

\noindent This definition extends to homology with coefficients and then to equivariant homology so that if $X=U \cup V$ and $M$ is a $\pi_1(X)$-module (hence also a $\pi_1(U)$-, $\pi_1(V)$-, and $\pi_1(U \cap V)$-module), there is a long exact sequence of equivariant homology groups with coefficients in $M$: $$\cdots \rightarrow H_n(U \cap V; M) \rightarrow H_n(U; M) \oplus H_n(V; M) \rightarrow H_n(X; M) \rightarrow\cdots$$

\section{Knot Polynomials}
Polynomials are an important category of knot invariants that can encode more subtle information about specific knots.  The Alexander Polynomial, discovered by J.W. Alexander in 1928, was the first of these polynomials.  While it can be computed directly from a presentation of a knot via a skein relation, I will focus, rather, on its homological formulation.  Unless otherwise directed, see \cite{Hillman}, \cite{Lickorish}, \cite{Murasugi} for basic references in this section.

\subsection{Elementary Ideals}

Many determinantal invariants of knots, and more generally of modules, are given in the form of an elementary ideal of some finite presentation of the module, defined as follows.

\begin{prop} \label{xxx}
Let $R$ be a commutative Noetherian ring and $M$ a finitely generated $R$-module.  Then an exact sequence: $$R^p \stackrel{Q}{\rightarrow} R^q \rightarrow M \rightarrow 0$$ is a finite presentation for $M$ with $p\times q$ presentation matrix $Q$.  The $r^{th}$ elementary ideal, $E_r(M)$ is the ideal of $R$ generated by all of the $(q-r) \times (q-r)$ minors of $Q$.  The smallest principal ideal of $M$ containing $E_r(M)$ is denoted $\widetilde{E}_r(M)$
\end{prop}

\noindent A natural question is whether a relationship between modules corresponds to a relationship between their determinantal invariants.  If it does, then perhaps certain properties are encoded in the invariants.  Theorem 3.12 of \cite{Hillman} establishes such a correspondence.

\newtheorem{thm}{Theorem}[subsection]
\begin{thm} \label{xxx}
Let $0\rightarrow K \rightarrow M \rightarrow C \rightarrow 0$ be an exact sequence of $R$-modules.  Then if $K$ is an $R$-torsion module: $$\widetilde{E}_r(M) = \widetilde{E}_0(K)\widetilde{E}_r(C)$$ where $r$ is the rank of $C$.
\end{thm}

\noindent For knots in particular, we will be concerned only with determinantal invariants of associated torsion modules.  A torsion module is a module over a ring for which every element of the module has a nonzero annihilator in the ring.  Thus, the rank of a torsion module is 0, and we have the following corollary.

\begin{cor} \label{3tor}
Let $0\rightarrow K \rightarrow M \rightarrow C \rightarrow 0$ be an exact sequence of $R$-torsion modules.  Then: $$\widetilde{E}_0(M) = \widetilde{E}_0(K)\widetilde{E}_0(C).$$
\end{cor}

\begin{cor} \label{4tor}
Let $0\rightarrow A \rightarrow B \stackrel{\alpha}{\rightarrow} C \stackrel{\beta}{\rightarrow} D \rightarrow 0$ be an exact sequence of $R$-torsion modules.  Then if $R$ is a domain: $$\widetilde{E}_0(A)\widetilde{E}_0(C) = \widetilde{E}_0(B)\widetilde{E}_0(D).$$
\end{cor}

\begin{proof} Note that: $$0\rightarrow A \rightarrow B \rightarrow Im(\alpha) \rightarrow 0$$ $$0\rightarrow Ker(\beta) \rightarrow C \rightarrow D \rightarrow 0$$ are exact sequences of $R$-torsion modules.  Then by Corollary \ref{3tor}: $$\widetilde{E}_0(B) = \widetilde{E}_0(A)\widetilde{E}_0(Im(\alpha))$$ $$\widetilde{E}_0(C) = \widetilde{E}_0(Ker(\beta))\widetilde{E}_0(D).$$  Cross multiplying: $$ \widetilde{E}_0(A)\widetilde{E}_0(C)\widetilde{E}_0(Im(\alpha)) = \widetilde{E}_0(B)\widetilde{E}_0(D)\widetilde{E}_0(Ker(\beta)).$$  Since $R$ is a domain and since $Im(\alpha)$ and $Ker(\beta)$ are torsion, $\widetilde{E}_0(Im(\alpha))$ and $\widetilde{E}_0(Ker(\beta))$ are nonzero.  Then, since $R$ is a domain and $Im(\alpha) = Ker(\beta)$, we can cancel to obtain: $$\widetilde{E}_0(A)\widetilde{E}_0(C) = \widetilde{E}_0(B)\widetilde{E}_0(D).$$
\end{proof}

\noindent Relationships between invariants of specific modules can be refined by explicit calculation of their elementary ideals.  The next section defines the necessary tools for these calculations.

\subsection{Free Differential Calculus}

The free differential calculus was invented and explored in a series of papers by Ralph Fox beginning in 1953 \cite{Fox}.  Fox's derivatives are functions on free groups that resemble, in certain ways, ordinary derivatives as in calculus.  They are defined axiomatically as follows.

\begin{prop} \label{xxx}
Let $G$ be a free group generated freely by $x_1, x_2, \ldots, x_n$. Then for any word $w$ in $G$, its free derivative in the free group ring $\Z G$ is computed using the following formulas:\\
\\
\indent $1)$ $\frac{\partial \textbf{1}}{\partial x_j} = 0$\\
\indent $2)$ $\frac{\partial x_k}{\partial x_j} = \delta_{j,k}$, the Kronecker delta\\
\indent $3)$ $\frac{\partial x_j^{-1}}{\partial x_j} = -x_j^{-1}$\\
\indent $4)$ $\frac{\partial uv}{\partial x_j} = \frac{\partial u}{\partial x_j} + u\frac{\partial v}{\partial x_j}$, for $u,v$ words in $G$\\
\indent $5)$ $w = \textbf{1} + \sum_j \frac{\partial w}{\partial x_j}(x_j - \textbf{1})$, the \textit{fundamental formula}\\
\\
\noindent Together, these axioms define mappings $\frac{\partial}{\partial x_j}: G\rightarrow \Z G$ that can be extended in an obvious way to mappings $\frac{\partial}{\partial x_j}: \Z G\rightarrow \Z G$.
\end{prop}

\noindent The next section will show that free derivatives arise naturally in the formulation of the Alexander polynomial of a knot.

\subsection{The Alexander Polynomial}
Let $K$ be a knot and $G = \pi_1(S^3 - K)$ its group.  Consider a Wirtinger presentation, $G = \left\langle x_1,\ldots, x_m | R_1,\ldots, R_{m-1}\right\rangle$, with deficiency 1 \cite{Murasugi}.  Following from the free differential calculus, there is a resolution of $\Z$ over $\Z G$: 

$$0 \longrightarrow \Z G^{m-1} \stackrel{\partial_2}{\longrightarrow} \Z G^{m} \stackrel{\partial_1}{\longrightarrow} \Z G \stackrel{x_i\mapsto 1}{\longrightarrow}\Z \longrightarrow 0$$

\noindent where $\partial_2$ is right multiplication by $A = \left(\frac{\partial R_i}{\partial x_j}\right)$, the Jacobian of Fox free derivatives, and $\partial_1$ is right multiplication by $(1-x_1,\ldots, 1-x_m)^T$.\\

\noindent Now let $\psi : G \rightarrow G/[G, G]$ be the abelianization map.  If we fix $x_1$ to be a meridian of $K$,  then $G/[G, G] \cong \left\langle x_1\right\rangle$, so we can assume $\psi x_j = 1$, $j\neq 1$.  Thus, there is a canonical $G$-action on $M = \Z [t^\pm]$, the ring of integer Laurent polynomials, under which $x_1$ acts on $M$ as multiplication by $t$, and $x_j$, $j\neq 1$ acts trivially.  Denote this action by $\Psi$.  Then tensoring by $M$ over $\Z G$ we have a complex: $$0 \rightarrow M^{m-1} \stackrel{\alpha_2}{\rightarrow} M^{m} \stackrel{\alpha_1}{\rightarrow} M \rightarrow 0$$

\noindent where $\alpha_2 = \Psi\partial_2$ is right multiplication by a matrix which we denote by $A^\Psi$ and $\alpha_1 = \Psi\partial_1$ is right multiplication by $(1-t,0,\ldots,0)^T$.

\begin{prop} \label{xxx}
The Alexander polynomial, $\Delta_K(t)$, of $K$ is a generator of $\widetilde{E}_0(M)$--the smallest principal ideal of $M$ containing $E_0(M)$, the ideal of $M$ containing all $(m-1) \times (m-1)$ minors of $A^\Psi$.
\end{prop}

\noindent Now we consider the homology of $X_K = S^3 - K$.

\begin{thm} \label{hap}
$\Delta_K(t)$ is a generator of $\widetilde{E}_0[H_1(X_K; M)]$ and is nonzero if and only if $H_2(X_K; M) = 0$.
\end{thm}

\begin{proof} From the definition of $\alpha_1$, we see that the first column of $A^\Psi$ (corresponding to $x_1$, the meridian of K) is 0.  Letting $\alpha_2'$ be given by right multiplication by $A^{\Psi '}$, the $(m-1) \times (m-1)$ matrix formed from $A^\Psi$ by removing the first column, we obtain the exact sequence:

$$0 \longrightarrow M^{m-1} \stackrel{\alpha_2'}{\longrightarrow} Ker(\alpha_1) \cong M^{m-1} \stackrel{\alpha_1}{\longrightarrow} \frac{Ker(\alpha_1)}{Im(\alpha_2)} = H_1(X_K; M) \longrightarrow 0$$

\noindent from which we deduce that $\widetilde{E}_0[H_1(X_K; M)] = det(A^{\Psi '}) = \widetilde{E}_0(M).$\\

\noindent We also deduce that $\Delta_K(t) \neq 0 \leftrightarrow det(A^{\Psi '}) \neq 0 \leftrightarrow Ker(\alpha_2') = 0 \leftrightarrow H_2(X_K; M) = 0$.
\end{proof}

\noindent We will also need to be able to compute the polynomial of a link $L$ of 2 components, $K$ and the unknot, $A$.  Proposition 2.3 of \cite{Murasugi2} states that a link $L$ also has a presentation $G_L = \left\langle x_1,\ldots, x_m | R_1,\ldots, R_{m-1}\right\rangle$ of deficiency 1.  Also, $G_L/[G_L, G_L] \cong \left\langle m_K, m_A\right\rangle$.  Thus, we can assume $x_1 = m_K$, $x_2 = m_A$, and $x_j \in [G_L, G_L]$ for $j>2$.  As explained above for $H_*(X_K; M)$, the groups $H_*(X_L; M)$ are computed from a complex: $$0\rightarrow M^{m-1}\rightarrow M^m \rightarrow M \rightarrow 0$$

\noindent and as before, $\Delta_L(t)$ is a generator of $\widetilde{E}_0[H_1(X_L; M)]$ and is nonzero if and only if $H_2(X_L; M) = 0$.

\subsection{The Twisted Alexander Polynomial}

The twisted Alexander polynomial for knots was discovered by X.S. Lin and generalized to finitely presentable groups by M. Wada.  Here, we will focus on Wada's formulation and then consider a more recent homological formulation.\\

\noindent As before, let $K$ be a knot and $G = \pi_1(S^3 - K)$ its group.  Consider a Wirtinger presentation, $G = \left\langle x_1,\ldots, x_m | R_1,\ldots, R_{m-1}\right\rangle$, with deficiency 1.  There is a resolution of $\Z$ over $\Z G$: 

$$0 \longrightarrow \Z G^{m-1} \stackrel{\partial_2}{\longrightarrow} \Z G^{m} \stackrel{\partial_1}{\longrightarrow} \Z G \stackrel{x_i\mapsto 1}{\longrightarrow}\Z \longrightarrow 0$$

\noindent where $\partial_2$ is right multiplication by $A = \left(\frac{\partial R_i}{\partial x_j}\right)$, the Jacobian of Fox free derivatives, and $\partial_1$ is right multiplication by $(1-x_1,\ldots, 1-x_m)^T$.\\

\noindent The abelianization of $G$ is, again, isomorphic to the group generated by the meridian of $K$, to which we fix $x_1$.  Now define a representation of $G$ by a map $G \rightarrow GL_n(R)$ which extends to a map $\rho: \Z G \rightarrow M_n(R)$, where $M_n(R)$ is the ring of matrices of order $n$ with entries in a U.F.D. $R$.  We now extend the canonical $G$-action of the usual Alexander polynomial to an action on $R[t^\pm]$ under which $x_1$ acts as multiplication by $\rho(x_1)t$, and $x_j$, $j \neq 1$ acts by $\rho(x_j)$.  Denote this action by $\Phi$.  Then tensoring by $R[t^\pm]$ over $\Z G$ we have:  

$$0 \longrightarrow R[t^\pm]^{n(m-1)} \stackrel{\beta_2}{\longrightarrow} R[t^\pm]^{nm} \stackrel{\beta_1}{\longrightarrow} R[t^\pm] \longrightarrow 0$$

\noindent where $\beta_2$ is right multiplication by the $n(m-1) \times nm$ matrix $A^\Phi$ and $\beta_1$ is right multiplication by the $nm \times n$ matrix $(I-\Phi x_1,\ldots, I-\Phi x_m)^T$.

\begin{prop} [Wada] \label{wada}
For some $j$, $det \Phi(1-x_j) \neq 0$.  Let $A_j$ be the $(m-1) \times (m-1)$ matrix obtained from $A$ by deleting the $j^{th}$ column (associated with the generator $x_j$).  Let $A_j^\Phi$ be the $n(m-1) \times n(m-1)$ matrix obtained from $A_j$ by the action $\Phi$.  The Wada twisted Alexander polynomial, $\Delta^W_{K,\rho}(t)$, of $K$ associated to the representation $\rho$ is given by the rational expression $\frac{det(A_j^\Phi)}{det \Phi(1-x_j)}$.  For a given knot and representation, it is invariant up to a unit factor in $R[t^\pm]$ \cite{Wada}.
\end{prop}

\noindent There are several different formulations of the twisted Alexander polynomial.  Though they don't all agree precisely, they are equivalent in the sense that they describe invariants of a given knot and presentation.  For example, Lin's original invariant corresponds to the numerator of Wada's.  More recently, a twisted invariant was described by Kirk and Livingston in terms of homology and related back to Definition \ref{wada} \cite{Kirk}.

\begin{prop} \label{tap}
The $i^{th}$ twisted Alexander polynomial of $K$ associated to $\rho$, denoted $\Delta^i_{K,\rho}(t)$, is a generator of $\widetilde{E}_0[H_i(X_K; R[t^\pm])]$.  For $i=1$ this is called the twisted Alexander polynomial and denoted $\Delta_{K,\rho}(t)$.  The invariant described by Wada (Definition \ref{wada}) is given by $\Delta^W = \frac{\Delta}{\Delta^0}$.
\end{prop}

\noindent Since there is also a deficiency 1 presentation for a link group, these definitions generalize to twisted link polynomials as in the previous section.

\section{Murasugi's Condition}

Originally published in 1971 by Kunio Murasugi, the following condition is one of the first unrestricted theorems on periodic knots.  It relates the Alexander Polynomial of such a knot to a cyclotomic polynomial:

\begin{thm} [Murasugi] \label{xxx}
If $K$ is a periodic knot of order $p^r$ in $S^3$, $p$ a prime, then the knot polynomial $\Delta_K(t)$ of $K$ must satisfy: $$\Delta_K(t) \equiv f(t)^{p^r}(1+ t + t^2 + \cdots t^{\lambda - 1})^{p^r-1}\hspace{5 pt} (mod~p)$$
for some knot polynomial $f(t)$ and a positive integer $\lambda$, $(\lambda, p) = 1$. 
\end{thm}

\noindent The topological interpretation of both $f(t)$ and $\lambda$ will be discussed in the following.  Murasugi's original proof relies heavily on specific presentations of knot groups and determining how they relate to the corresponding presentation matrices via the free differential calculus.  This sort of argument proves unwieldy for extension to the twisted polynomial.  However, we have seen that the Alexander polynomial of a knot can be expressed solely in terms of the homology of the universal cover of its complement in $S^3$.  Thus, I provide a new proof using homology theory which will be sufficiently abstract to allow generalization of Murasugi's condition to the twisted case. 

\subsection{A Homological Proof}

\subsubsection{Periodic Knots and Cyclic Covers}

\begin{prop} \label{xxx}
Let $\Sigma$ be a homology 3-sphere.  A knot $K$ in $\Sigma$ has period $q$, or is \textit{periodic} of order $q$, if there is an orientation-preserving homeomorphism $\phi: \Sigma\rightarrow\Sigma$ under which:\\
\\
\indent $1)$ The set of fixed points is a circle (the unknot) $A$ disjoint from $K$;\\
\indent $2)$ $\phi(K) = K$;\\
\indent $3)$ $\phi^q = \textbf{1}$ but $\phi^{q'} \neq \textbf{1}$ for $0<q'<q$.\\
\end{prop}

\noindent Consider the orbit space $\overline{\Sigma} = \Sigma/\phi$ which is also a homology 3-sphere \cite{Murasugi2}.  The space $\Sigma$ together with the quotient map $\varrho: \Sigma\rightarrow\overline{\Sigma}$ is a $q$-fold cyclic cover of $\overline{\Sigma}$ branched along $\bar{A} = \varrho(A)$.  Let $\bar{K} = \varrho(K)$, $\bar{L} = \bar{K} \sqcup \bar{A}$, and $L = K \sqcup A$.  Furthermore, let $G_L = \pi_1(X_L)$ where $X_L = \Sigma - L$ and $G_{\bar{L}} = \pi_1(X_{\bar{L}})$.  Then the quotient group $G_{\bar{L}}/G_L \cong \Z/q\Z := \Z_q$ is generated by the meridian of $\bar{A}.$\\

\noindent Let $q=p^r$, $p$ prime, and consider $H_n(X_L; M_p)$ where: $$M_p = M/pM = \Z_p[t^\pm]$$ 

\noindent so that the coefficients of the usual $\Z[t^\pm]$ have been passed to $\Z_p:=\Z/p\Z$.  Noting that $M_p$ is a $G_{\bar{L}}$-module (and therefore a $G_L$-module), we apply Shapiro's Lemma to obtain:

$$H_n(X_L; M_p) \cong H_n(X_{\bar{L}}; M_p\otimes_{\Z_p} \Z_p[G_{\bar{L}}/G_L]).$$

\subsubsection{Homology and Ideals}

Now note that:

$$\Z_p[G_{\bar{L}}/G_L]\cong\Z_p[\Z_q]\cong\frac{\Z_p[x]}{(x^q-1)} = \frac{\Z_p[x]}{(x-1)^q}$$

\noindent where the isomorphism to the polynomial ring is given by mapping $g$, a generator of $G_{\bar{L}}/G_L$, to $x$.  The final equality holds since the coefficient ring is $\Z_p$.  For $k>1$, modules of the form $\frac{\Z_p[x]}{(x-1)^k}$ fit into a short exact sequence:

$$0 \rightarrow \frac{\Z_p[x]}{(x-1)^{k-1}} \stackrel{}{\rightarrow} \frac{\Z_p[x]}{(x-1)^k} \rightarrow \frac{\Z_p[x]}{(x-1)}\cong\Z_p \rightarrow 0$$

\noindent where the second map is multiplication by $1-x$.

\begin{lem} \label{rec}
Let $N_k = \frac{\Z_p[x]}{(x-1)^k}$.  Then: $$\widetilde{E}_0[H_1(X_{\bar{L}}; M_p \otimes N_k)]\cdot \widetilde{E}_0[H_0(X_{\bar{L}}; M_p)] = \widetilde{E}_0[H_1(X_{\bar{L}}; M_p \otimes N_{k-1})] \cdot \widetilde{E}_0[H_1(X_{\bar{L}}; M_p)]$$
\end{lem}

\begin{proof} Noting that $M_p \otimes N_1 \cong M_p$, we have the following long exact sequence:

\begin{flushleft}
$\cdots \rightarrow H_2(X_{\bar{L}}; M_p) \rightarrow H_1(X_{\bar{L}}; M_p \otimes N_{k-1}) \rightarrow H_1(X_{\bar{L}}; M_p \otimes N_k)$
\end{flushleft}
\begin{flushright}
$\rightarrow H_1(X_{\bar{L}}; M_p) \rightarrow H_0(X_{\bar{L}}; M_p \otimes N_{k-1}) \rightarrow \cdots \rightarrow 0$
\end{flushright} \mbox{}\\

\noindent We can reduce this sequence using the following facts:

\noindent \textbf{1) $H_2(X_{\bar{L}}; M_p) = 0$}
\begin{tabbing}
\hspace{20 pt}
\= By a classical result of Seifert, $\Delta_K(0) = \pm 1$ for a knot $K$, so $\Delta_K(t) \neq 0$.\\
\> Then we have that $H_2(X_K; M) = 0$, so certainly $H_2(X_K; M_p) = 0$.\\
\> Also note that for $n > k$ the $n^{th}$ homology of a $k$-dimensional space is 0,\\ 
\> so $H_2(A; M_p) = H_3(X_K; M_p) = 0$.\\
\> \\
\> Consider the decomposition $X_K = X_L \cup N(A)$ where $N(A)$ is an open\\ 
\> tubular neighborhood of $A$.  Since $A$ is a homotopy retract of $N(A)$ and $\T^2$,\\
\> the 2-torus, is a homotopy retract of $N(A) \cap X_L$, $H_n(N(A); M_p) =$\\
\> $H_n(A; M_p)$ and $H_n(N(A) \cap X_L; M_p) = H_n(\T^2; M_p)$.  Therefore, the\\  
\> Mayer-Vietoris sequence becomes:\\
\> \\
\hspace{28 pt} $0\rightarrow H_2(\T^2; M_p) \rightarrow H_2(X_L; M_p) \oplus H_2(A; M_p)=0 \rightarrow H_2(X_K; M_p)=0$\\
\\
\> and then:\\
\>\\
\hspace{97 pt} $0\rightarrow H_2(\T^2; M_p) \rightarrow H_2(X_L; M_p) \rightarrow 0$\\
\> \\
\> So $H_2(X_L; M_p) = H_2(\T^2; M_p)$.  The fundamental group, $\pi_1(\T^2) \cong \Z \oplus \Z$, is\\
\> generated by $m_A$ and $l_A$, the meridian and longitude of A.  Furthermore the\\
\> torus $\T^2$ is a CW-complex with 1 0-cell, 2 1-cells, and 1 2-cell.  Thus, the\\
\> homology groups $H_*(\T^2; M_p)$ can be computed by tensoring the sequence:\\
\\ 
\hspace{80 pt} $0\rightarrow \Z[\Z\oplus\Z]\stackrel{\gamma_2}{\rightarrow}\Z[\Z\oplus\Z]^2\stackrel{\gamma_1}{\rightarrow}\Z[\Z\oplus\Z]\rightarrow 0$\\
\\
\> with $M_p$ over $\Z[\Z\oplus\Z] = \Z[\pi_1(\pi^2)]$.  Here $\gamma_2$ is multiplication by\\
\> $(1-l_A, m_A-1)$ and $\gamma_1$ is multiplication by $(1-m_A, 1-l_A)^T$.  These maps\\
\> come from the free differential calculus and the presentation\\
\> $\pi_1(\pi^2) \cong \left\langle m_A, l_A| m_Al_Am_A^{-1}l_A^{-1} \right\rangle$.  Therefore, by direct computation,\\
\> $H_2(\T^2; M_p)$ is 0, so $H_2(X_L; M_p) = 0$\\
\> \\
\> Applying Shapiro's Lemma, $H_2(X_{\bar{L}}; M_p \otimes N_q) = 0$.\\
\> \\
\> Noting again that in our long exact sequence the $H_3$ groups are 0, then\\
\> $H_2(X_{\bar{L}}; M_p \otimes N_{k-1}) \rightarrow H_2(X_{\bar{L}}; M_p \otimes N_k)$ is an injection.  Thus, $H_2(X_{\bar{L}}; M_p)$\\ 
\> injects into $H_2(X_{\bar{L}}; M_p \otimes N_q) = 0$, and the desired result is obtained.
\end{tabbing}

\noindent\textbf{2) $H_0(X_{\bar{L}}; M_p \otimes N_k) \rightarrow H_0(X_{\bar{L}}; M_p)$} is an isomorphism
\begin{tabbing}
\hspace{20 pt}
\= First note that this map is a surjection since the following map in the se-\\
\> quence is trivial.  Now we must show it is injective.  As explained in Section\\
\> 3.3, the groups $H_*(X_{\bar{L}}; M_p \otimes N_k)$ are computed from a complex:\\
\\
\hspace{67 pt} $0\rightarrow (M_p \otimes N_{k})^{m-1} \rightarrow (M_p \otimes N_{k})^m \rightarrow M_p \otimes N_{k} \rightarrow 0$\\
\\
\> where if $x_1,\ldots,x_m$ are generators of $G_{\bar{L}}$, the second to last map is right\\
\> multiplication by $(1-x_1,\ldots,1-x_m)^T$.  By the prescribed $G_{\bar{L}}$-actions on\\
\> both $M_p$ and $N_k$, $m_K$ acts by multiplication by $t \otimes 1$ and $m_A$ acts by mult-\\
\> iplication by $1 \otimes x$.  The remaining $x_i$ act trivially.\\
\\
\> Thus, $(1-(1\otimes x))(M_p \otimes N_k) = M_p\otimes N_k - M_p\otimes xN_k = M_p\otimes (1-x)N_k$\\
\> has trivial image in $H_0(X_{\bar{L}}; M_p \otimes N_k)$.  And since the map:\\
\\
\hspace{95 pt} $H_0(X_{\bar{L}}; M_p \otimes N_{k-1}) \rightarrow H_0(X_{\bar{L}}; M_p \otimes N_{k})$\\
\\
\> is given by multiplication by $1-x$, it must be the zero map.  So the map\\
\> $H_0(X_{\bar{L}}; M_p \otimes N_k) \rightarrow H_0(X_{\bar{L}}; M_p)$ is injective and therefore an isomorphism.
\end{tabbing}

\noindent These two facts reduce our long exact sequence to: $$0 \rightarrow H_1(X_{\bar{L}}; M_p \otimes N_{k-1}) \rightarrow H_1(X_{\bar{L}}; M_p \otimes N_k) \rightarrow H_1(X_{\bar{L}}; M_p) \rightarrow H_0(X_{\bar{L}}; M_p) \rightarrow 0$$

\noindent Then by Corollary \ref{4tor}: $$\widetilde{E}_0[H_1(X_{\bar{L}}; M_p \otimes N_k)] \widetilde{E}_0[H_0(X_{\bar{L}}; M_p)] = \widetilde{E}_0[H_1(X_{\bar{L}}; M_p \otimes N_{k-1})] \widetilde{E}_0[H_1(X_{\bar{L}}; M_p)].$$
\end{proof}

\noindent Lemma \ref{rec} gives a recursion formula which can be used in a downward inductive argument to produce the following corollary.

\begin{cor} \label{ind}
$\widetilde{E}_0[H_1(X_{L}; M_p)](\widetilde{E}_0[H_0(X_{\bar{L}}; M_p)])^{q-1}  = (\widetilde{E}_0[H_1(X_{\bar{L}}; M_p)])^q$
\end{cor}

\begin{proof} By Lemma \ref{rec} we have:
$$\widetilde{E}_0[H_1(X_{\bar{L}}; M_p \otimes N_k)] \widetilde{E}_0[H_0(X_{\bar{L}}; M_p)] = \widetilde{E}_0[H_1(X_{\bar{L}}; M_p \otimes N_{k-1})] \widetilde{E}_0[H_1(X_{\bar{L}}; M_p)]$$
Shifting $k$ to $k-1$ (and therefore $k-1$ to $k-2$), a similar formula is obtained.  Combining the two:
$$\widetilde{E}_0[H_1(X_{\bar{L}}; M_p \otimes N_k)](\widetilde{E}_0[H_0(X_{\bar{L}}; M_p)])^2 = \widetilde{E}_0[H_1(X_{\bar{L}}; M_p \otimes N_{k-2})] (\widetilde{E}_0[H_1(X_{\bar{L}}; M_p)])^2$$
Thus, in general, we can continue this downward inductive process on $k$ to obtain:
$$\widetilde{E}_0[H_1(X_{\bar{L}}; M_p \otimes N_k)] (\widetilde{E}_0[H_0(X_{\bar{L}}; M_p)])^l = \widetilde{E}_0[H_1(X_{\bar{L}}; M_p \otimes N_{k-l})]  (\widetilde{E}_0[H_1(X_{\bar{L}}; M_p)])^l$$
Letting $k = q$, $l = q-1$ and applying Shapiro's Lemma on the first term:
$$\widetilde{E}_0[H_1(X_{L}; M_p)](\widetilde{E}_0[H_0(X_{\bar{L}}; M_p)])^{q-1}  = (\widetilde{E}_0[H_1(X_{\bar{L}}; M_p)])^q$$
\end{proof}

\subsubsection{Relating $\Delta_K(t)$ to $\Delta_L(t)$}

\begin{lem} \label{rel}
$\Delta_L (t) = (1-t^\lambda)\Delta_K(t)$ and $\Delta_{\bar{L}} (t) = (1-t^\lambda)\Delta_{\bar{K}}(t)$ where $\lambda = lk(K, A)$.
\end{lem}

\begin{proof} Consider the Mayer-Vietoris sequence for $X_K = X_L \cup N(A)$ with coefficients in $M = \Z[t^\pm]$: 

\begin{flushleft}
$0\rightarrow H_1(\T^2; M) \rightarrow H_1(X_L; M) \oplus H_1(A; M) \rightarrow H_1(X_K; M)$
\end{flushleft}
\begin{flushright}
$\rightarrow H_0(\T^2; M) \rightarrow H_0(X_L; M) \oplus H_0(A; M) \rightarrow H_0(X_K; M) \rightarrow 0$
\end{flushright}

\noindent Since the fundamental group of $A$ is the free cyclic group generated by $l_A$ (the longitude of $A$), a chain complex computing the equivariant homology of $A$ is: $$0 \rightarrow \Z \pi_1(A) \stackrel{1-l_A}{\longrightarrow} \Z \pi_1(A) \rightarrow 0$$

\noindent Notice that $l_A$ traverses the meridian of $K$ $\lambda = lk(K, A)$ times.  This is true since the linking number, by definition, represents the minimum number of times two knots need to pass through each other to separate, or the number of times they wind around each other.  Thus, $l_A = m_K^\lambda$ in $G_K$, and so the action of $\pi_1(A)$ on $M$ is such that $l_A$ acts as multiplication by $t^\lambda$.  Tensoring the above sequence with $M$: $$0 \rightarrow M \stackrel{1-t^\lambda}{\longrightarrow} M \rightarrow 0$$

\noindent is the presentation complex for $A$ which corresponds to its equivariant homology.  Since $\lambda$ is nonzero the map given by $1-t^\lambda$ has no kernel, so $H_1(A; M) = 0$ and $H_0(A; M) = M/(1-t^\lambda)M$.  As explained in the proof of Lemma \ref{rec}, the homology groups $H_*(\T^2; M)$ are computed by a complex: $$0\rightarrow M \stackrel{\gamma_2}{\rightarrow} M^2 \stackrel{\gamma_1}{\rightarrow} M\rightarrow 0$$

\noindent where $\gamma_2$ is right multiplication by $(1-t^\lambda, 0)$ and $\gamma_1$ is right multiplication by $(0, 1-t^\lambda)^T$.  Thus, $H_0(\T^2; M) \rightarrow H_0(A; M)$ is an isomorphism, and $H_1(\T^2; M) \cong \frac{M}{(1-t^\lambda)M}$.  In particular, $\widetilde{E}_0[H_1(\T^2; M)] = (1-t^\lambda)$. The Mayer-Vietoris sequence is now reduced to: $$0 \rightarrow H_1(\T^2; M) \rightarrow H_1(X_L; M) \rightarrow H_1(X_K; M) \rightarrow 0.$$

\noindent So by Corollary \ref{3tor}: $$\widetilde{E}_0[H_1(X_L; M)] = \widetilde{E}_0[H_1(\T^2; M)] \cdot \widetilde{E}_0[H_1(X_K; M)].$$

\noindent Thus, we obtain: $$\Delta_L (t) = (1-t^\lambda)\Delta_K(t).$$

\noindent By the same argument, this relationship also holds for $\Delta_{\bar{L}}$ and $\Delta_{\bar{K}}$.  To see that the linking numbers $lk(K, A)$ and $lk(\bar{K}, \bar{A})$ are the same, note that the covering map $\varrho$ induces a homomorphism $\varpi: G_L\rightarrow G_{\bar{L}}$.  Using $l_A = m_K^\lambda$ and homomorphism properties, $\varpi l_A = \varpi(m_K^\lambda) = (\varpi m_K)^\lambda$.  Since under $\varpi$, $l_A\mapsto l_{\bar{A}}$ and $m_K\mapsto m_{\bar{K}}$, then $l_{\bar{A}} = m_{\bar{K}}^\lambda$ and the definition of linking number implies that $\lambda = lk(\bar{K}, \bar{A})$.
\end{proof}

\subsubsection{Murasugi's Condition}

\begin{thm} [Murasugi] \label{xxx}
$\Delta_K(t) \equiv \Delta_{\bar{K}}(t)^q \left(\frac{1-t^\lambda}{1-t}\right)^{q-1}\hspace{5 pt} (mod~p)$
\end{thm}

\begin{proof} By Corollary \ref{ind} we have:
$$\widetilde{E}_0[H_1(X_{L}; M_p)](\widetilde{E}_0[H_0(X_{\bar{L}}; M_p)])^{q-1}  = (\widetilde{E}_0[H_1(X_{\bar{L}}; M_p)])^q$$

\noindent If $M=\Z[t^\pm]$, then $M_p=M/p$ and a presentation complex for $M_p$ is given by a short exact sequence: $$0\rightarrow M\stackrel{\times p}{\rightarrow} M \rightarrow M_p \rightarrow 0.$$

\noindent This corresponds to a long exact sequence of homology groups for $X_{\bar{L}}$, which (since $H_2(X_{\bar{L}}; M_p) = 0$) reduces to:

\begin{flushleft}
\hspace{15 pt} $0\rightarrow H_1(X_{\bar{L}}; M) \stackrel{\times p}{\rightarrow} H_1(X_{\bar{L}}; M) \rightarrow H_1(X_{\bar{L}}; M_p)$
\end{flushleft}
\begin{flushright}
$\rightarrow H_0(X_{\bar{L}}; M)\stackrel{\times p}{\rightarrow} H_0(X_{\bar{L}}; M) \rightarrow H_0(X_{\bar{L}}; M_p)\rightarrow 0.\hspace{15 pt}$
\end{flushright}

\noindent From the presentation $M^{m-1} \rightarrow M^m  \stackrel{\times(1-t)}{\longrightarrow} M \rightarrow 0$, the group $H_0(X_{\bar{L}}; M)$ is computable as $\frac{M}{(1-t)M}$.  Notice that $H_0(X_{\bar{L}}; M)\stackrel{\times p}{\rightarrow} H_0(X_{\bar{L}}; M)$ is an inclusion $\frac{M}{(1-t)M} \stackrel{\times p}{\hookrightarrow} \frac{M}{(1-t)M}$. So the cokernel of the map $H_1(X_{\bar{L}}; M_p)\rightarrow H_0(X_{\bar{L}}; M)$ is 0, and this long exact sequence can be broken into the short exact sequences: $$0\rightarrow H_1(X_{\bar{L}}; M) \stackrel{\times p}{\rightarrow} H_1(X_{\bar{L}}; M) \rightarrow H_1(X_{\bar{L}}; M_p)\rightarrow 0$$ $$0\rightarrow H_0(X_{\bar{L}}; M) \stackrel{\times p}{\rightarrow} H_0(X_{\bar{L}}; M) \rightarrow H_0(X_{\bar{L}}; M_p)\rightarrow 0$$

\noindent which imply $H_1(X_{\bar{L}}; M_p)\cong H_1(X_{\bar{L}}; M)/p$ and $H_0(X_{\bar{L}}; M_p)\cong H_0(X_{\bar{L}}; M)/p$.//

\noindent This trivially extends to: $$\widetilde{E}_0[H_1(X; M_p)] \equiv \widetilde{E}_0[H_1(X; M)]\hspace{5 pt} (mod~p)$$ $$\widetilde{E}_0[H_0(X; M_p)] \equiv \widetilde{E}_0[H_0(X; M)]\hspace{5 pt} (mod~p).$$

\noindent Thus, we have: $$\widetilde{E}_0[H_1(X_{L}; M)](\widetilde{E}_0[H_0(X_{\bar{L}}; M)])^{q-1}  \equiv (\widetilde{E}_0[H_1(X_{\bar{L}}; M)])^q \hspace{5 pt} (mod~p).$$

\noindent Noting that products of ideals are generated by the products of their generators, Theorem \ref{hap} gives us: $$\Delta_L(t)(\widetilde{E}_0[H_0(X_{\bar{L}}; M)])^{q-1}  \equiv \Delta_{\bar{L}}(t)^q\hspace{5 pt} (mod~p).$$

\noindent Finally, by direct computation, $\widetilde{E}_0[H_0(X_{\bar{L}}; M)] = 1-t$.  Thus, by Lemma \ref{rel}: $$\Delta_K(t)(1-t^\lambda)(1-t)^{q-1} \equiv (\Delta_{\bar{K}}(t)(1-t^\lambda))^q \hspace{5 pt} (mod~p).$$

\noindent Hence: $$\Delta_K(t) \equiv \Delta_{\bar{K}}(t)^q \left(\frac{1-t^\lambda}{1-t}\right)^{q-1}\hspace{5 pt} (mod~p).$$

\noindent Thus, Murasugi's condition holds with $f(t) = \Delta_{\bar{K}}(t).$
\end{proof}

\subsection{The Twisted Case}

As noted, the twisted Alexander polynomial is not only an invariant of a knot but also of the choice of representation for its group.  Noting that the homological proof of Murasugi's condition involves calculations using mainly $mod~p$ coefficients, we restrict our consideration to $mod~p$ representations and see that an extended condition follows rather easily.\\

\noindent Keeping the notation from Section 4.1, let $\rho: G_K \rightarrow GL_n(\Z_p)$ be a representation for $G_K$ and $\bar{\rho}: G_{\bar{K}} \rightarrow GL_n(\Z_p)$ be the associated representation for $G_{\bar{K}}$, the group of the quotient knot.\\

\noindent The twisted homology groups will then be computed with coefficients in $R[t^\pm] = \Z_p[t^\pm] = M_p$, and we obtain the following results.

\begin{lem} \label{trel}
$\Delta_{L, \rho}(t) = \det(I_n - \rho(l_A)t^\lambda)\Delta_{K, \rho}(t)$ and $\Delta_{\bar{L}, \bar{\rho}}(t) = \det(I_n - \rho(l_A)t^\lambda)\Delta_{\bar{K}, \bar{\rho}}(t)$ where $\lambda = lk(K, A)$.  Furthermore, $\Delta_{\bar{L}, \bar{\rho}}^0 \cong \Delta_{\bar{K},\rho}^0$.
\end{lem}

\begin{proof}
Consider the Mayer-Vietoris sequence for $X_K = X_L \cup N(A)$ with coefficients in $M_p = \Z_p[t^\pm]$: 

\begin{flushleft}
$0\rightarrow H_1(\T^2; M_p) \rightarrow H_1(X_L; M_p) \oplus H_1(A; M_p) \rightarrow H_1(X_K; M_p)$
\end{flushleft}
\begin{flushright}
$\rightarrow H_0(\T^2; M_p) \rightarrow H_0(X_L; M_p) \oplus H_0(A; M_p) \rightarrow H_0(X_K; M_p) \rightarrow 0$
\end{flushright}

\noindent The chain complex for $A$ is given by: $$0 \rightarrow \Z \pi_1(A) \stackrel{1-l_A}{\longrightarrow} \Z \pi_1(A) \rightarrow 0$$

\noindent where $l_A$ is the longitude of the unknot, which freely generates $G_A$.  As before, $l_A$ traverses the meridian of $K$ $\lambda = lk(K, A)$ times and thus, $l_A = m_K^\lambda$.  However, the action of $l_A$ on $M_p$ is now given by $\rho(l_A)t^\lambda$, so tensoring with $M_p$: $$0 \longrightarrow M_p \stackrel{I_n-\rho(l_A)t^\lambda}{\longrightarrow} M_p \longrightarrow 0$$

\noindent is the presentation complex for $A$ which corresponds to its equivariant homology with coefficients in $M_p$.  The homology groups $H_*(\T^2; M_p)$ are computed by a complex: $$0\rightarrow M_p \stackrel{\gamma_2}{\rightarrow} M_p^2 \stackrel{\gamma_1}{\rightarrow} M_p\rightarrow 0$$

\noindent where $\gamma_2$ is right multiplication by $(I_n-\rho(l_A)t^\lambda, 0)$ and $\gamma_1$ is right multiplication by $(0, I_n -\rho(l_A)t^\lambda)^T$.  Thus, $H_0(\T^2; M_p) \rightarrow H_0(A; M_p)$ is an isomorphism, and $H_1(\T^2; M_p) \cong \frac{M_p}{(I_n-\rho(l_A)t^\lambda)M_p}$.  In particular, $\widetilde{E}_0[H_1(\T^2; M_p)] = det(I_n-\rho(l_A)t^\lambda)$.  The Mayer-Vietoris sequence is now reduced to: $$0 \rightarrow H_1(\T^2; M_p) \rightarrow H_1(X_L; M_p) \rightarrow H_1(X_K; M_p) \rightarrow 0$$

\noindent So by Corollary \ref{3tor}: $$\widetilde{E}_0[H_1(X_L; M_p)] = \widetilde{E}_0[H_1(\T^2; M_p)] \widetilde{E}_0[H_1(X_K; M_p)]$$

\noindent Thus, we obtain: $$\Delta_{L, \rho}(t) = \det(I_n - \rho(l_A)t^\lambda)\Delta_{K, \rho}(t)$$

\noindent By the same argument, $\Delta_{\bar{L}, \bar{\rho}}(t) = \det(I_n - \bar{\rho}(l_{\bar{A}})t^{\bar{\lambda}})\Delta_{\bar{K}, \bar{\rho}}(t)$.  We have already shown that the linking numbers $\lambda$ and $\bar{\lambda}$ are the same.  Noting that $\rho(G_K) \subset \bar{\rho}(G_{\bar{K}})$ and that the quotent map is branched along $\bar{A}$, we obtain $\bar{\rho}(l_{\bar{A}}) = \rho(l_A)$.  Thus, $\Delta_{\bar{L}, \bar{\rho}}(t) = \det(I_n - \rho(l_A)t^\lambda)\Delta_{\bar{K}, \bar{\rho}}(t)$.\\

\noindent Finally, if we had instead considered the Mayer-Vietoris sequence for $X_{\bar{K}} = X_{\bar{L}} \cup N(\bar{A})$, we would have obtained the isomorphism $H_0(\T^2; M_p) \rightarrow H_0(\bar{A}; M_p)$ as before.  This time eliminating the $H_1$ groups and working with the $H_0$ groups, we see that $H_0(X_{\bar{L}}; M_p) \cong H_0(X_{\bar{K}}; M_p)$, so $\Delta_{\bar{L}, \bar{\rho}}^0 \cong \Delta_{\bar{K},\rho}^0$.
\end{proof}

\begin{thm} \label{xxx}
$\Delta_{K, \rho}(t) = \Delta_{\bar{K}, \bar{\rho}}(t)^q\left(\frac{det(I_n - \rho(l_A)t^\lambda)}{\Delta_{\bar{K},\bar{\rho}}^0}\right)^{q-1}$.\\  Alternatively, $\Delta_{K,\rho}(t) = \Delta_{\bar{K}, \bar{\rho}}(t) \left(\Delta^W_{\bar{K}, \bar{\rho}}(t) \det(I_n - \rho(l_A)t^\lambda)\right)^{q-1}.$ 
\end{thm}

\begin{proof}
In reviewing the arguments in the last section leading up to Corollary \ref{ind}, only the assumption that $H_2(X_K; M_p) = 0$ was related specifically to the Alexander polynomial. By Theorem \ref{hap}, this assumption did not sacrifice the generality of the argument.  In fact, this will also hold for the twisted polynomial since $H_2(X_K; M_p)$ is a free $M_p$-submodule of $C_2(X_K; M_p)$ with the same rank as $H_1(X_K; M_p)$.  Noticing that $H_1(X_K; M_p) = H_1(X_K; M) \otimes \Z_p$, which is torsion (since $\Delta_K(1) = \pm1$ so $\Delta_K$ $mod$ $p$ is nonzero), we have our result (thanks to Stefan Friedl for pointing this out).\\

\noindent Having resolved this issue, we can apply Corollary \ref{ind} to obtain: $$\widetilde{E}_0[H_1(X_{L}; M_p)](\widetilde{E}_0[H_0(X_{\bar{L}}; M_p)])^{q-1}  = (\widetilde{E}_0[H_1(X_{\bar{L}}; M_p)])^q$$

\noindent By Definition \ref{tap} this is: $$\Delta_{L, \rho}(t)(\Delta_{\bar{L}, \bar{\rho}}^0(t))^{q-1} = (\Delta_{\bar{L}, \bar{\rho}}(t))^q.$$

\noindent Then by a simple application of Lemma \ref{trel}: $$\Delta_{K, \rho}(t)\det(I_n - \rho(l_A)t^\lambda)(\Delta_{\bar{K}, \bar{\rho}}^0(t))^{q-1} = (\Delta_{\bar{K}, \bar{\rho}}(t)\det(I_n - \rho(l_A)t^\lambda))^q.$$

\noindent Rearranging: $$\Delta_{K,\rho}(t) = \Delta_{\bar{K}, \bar{\rho}}(t)^q\left(\frac{\det(I_n - \rho(l_A)t^\lambda)}{\Delta_{\bar{K}, \bar{\rho}}^0(t)}\right)^{q-1}.$$

\noindent Alternatively, we recall from Definition \ref{tap} that Wada's invariant satisies $\Delta^W = \frac{\Delta}{\Delta^0}$, so: $$\Delta_{K,\rho}(t) = \Delta_{\bar{K}, \bar{\rho}}(t) \left(\Delta^W_{\bar{K}, \bar{\rho}}(t) \det(I_n - \rho(l_A)t^\lambda)\right)^{q-1}.$$

\noindent Thus, the twisted extension of Murasugi's condition holds with $f(t) = \Delta_{\bar{K}, \bar{\rho}}(t)$. 

\end{proof}

\end{document}